# SPECTRAL PROPERTIES OF THE TANDEM JACKSON NETWORK, SEEN AS A QUASI-BIRTH-AND-DEATH PROCESS

By D. P. Kroese[1], W. R. W. Scheinhardt[2] and P. G. Taylor[1]

*University of Queensland, University of Twente and Centre for Mathematics and Computer Science, and University of Melbourne*

Quasi-birth-and-death (QBD) processes with infinite "phase spaces" can exhibit unusual and interesting behavior. One of the simplest examples of such a process is the two-node tandem Jackson network, with the "phase" giving the state of the first queue and the "level" giving the state of the second queue.

In this paper, we undertake an extensive analysis of the properties of this QBD. In particular, we investigate the spectral properties of Neuts's $R$-matrix and show that the decay rate of the stationary distribution of the "level" process is not always equal to the convergence norm of $R$. In fact, we show that we can obtain any decay rate from a certain range by controlling only the transition structure at level zero, which is independent of $R$.

We also consider the sequence of tandem queues that is constructed by restricting the waiting room of the first queue to some finite capacity, and then allowing this capacity to increase to infinity. We show that the decay rates for the finite truncations converge to a value, which is not necessarily the decay rate in the infinite waiting room case.

Finally, we show that the probability that the process hits level $n$ before level 0 given that it starts in level 1 decays at a rate which is not necessarily the same as the decay rate for the stationary distribution.

**1. Introduction.** A quasi-birth-and-death (QBD) process is a two-dimensional continuous-time Markov chain for which the generator has a block-tridiagonal structure. The first component of the QBD process is called the *level*, the second component the *phase*.

Received January 2003; revised October 2003.

[1]Supported by the Australian Research Council through Discovery Grant DP0209921.

[2]Supported by the Thomas Stieltjes Institute for Mathematics and the Netherlands Organisation for Scientific Research.

*AMS 2000 subject classification.* 60J27.

*Key words and phrases.* Decay rate, tandem Jackson network, QBD process, stationary distribution, hitting probabilities.

A comprehensive discussion of the properties of QBD processes with finitely many possible values of the phase variable can be found in the monographs of Neuts [11] and Latouche and Ramaswami [7]. In particular, it is known that the level process of a positive-recurrent QBD process with a finite phase space possesses a stationary distribution which decays geometrically as the level is increased. The decay parameter is equal to the spectral radius of Neuts's $R$-matrix, which is strictly less than 1. Similarly, the probability that a QBD process hits level $n$ before level 0 given that it starts in level 1 is known to decay geometrically with the same parameter.

For QBD processes with an *infinite* phase space the situation becomes more complicated. The $R$-matrix is now infinite-dimensional, and its spectral properties are not obvious. Also, the relationship between various decay parameters is different from the finite-dimensional case. A start in the study of such processes was made by Takahashi, Fujimoto and Makimoto [17]. They gave conditions under which the infinite-dimensional $R$ is $\alpha$-positive (see Section 2 for a definition). Under a further condition on the stationary distribution at level 0, they were then able to infer that the stationary distribution decays at rate $\alpha$. However, as we shall see in Section 4, there are many circumstances where the conditions of [17] are not satisfied.

The purpose of this paper is to make a contribution to the study of the behavior of infinite-phase QBD processes by considering a special case which exhibits interesting behavior. This special case is a two-node tandem Jackson network, in which the number of customers in the first queue gives the phase variable and the number of customers in the second queue gives the level variable. This system was studied via simulation in [6], where the authors used some of the results of the current paper to calculate the relevant decay rates. It is also a special case of the system studied in [3].

We show that, when the first queue has an infinite waiting room, the decay rate of the stationary distribution of the "level" process (the state of the second queue) may not be equal to the convergence norm of $R$, which can be thought of as the analogue of the spectral radius in the infinite-dimensional case. In fact, we show that we can construct a range of decay rates for the stationary distribution of the second queue by controlling only the transition structure when the second queue is empty, that is, at level 0. Futhermore, the decay rate, as $n \to \infty$, of the probability that the number of customers in the second queue hits $n$ before 0 given that it starts at 1 may not be the same as the decay rate of the stationary distribution. Such behavior does not occur in finite-phase QBD processes.

We also consider the limiting behavior of the tandem queue when the waiting room of the first queue is finite, and increases to infinity. We show that the eigenvalues of the $R$-matrix converge to a continuum, possibly with one additional isolated point—the latter being the case when the second buffer is the bottleneck. A consequence of this is that the decay rate in the infinite waiting room case may not be the same as the limiting value of the decay rates in the finite waiting room case.

The rest of the paper is organized as follows. In Section 2 we present some general results for QBD processes. We consider processes with both finite and infinite phase spaces. In Section 3 we formulate the two-node tandem Jackson network as a QBD process. In Sections 4 and 5, we discuss the decay rate of the stationary distribution when the capacity of the first queue is infinite and finite, respectively. These sections make heavy use of the properties of certain orthogonal polynomials. In Section 6, we show how we can obtain any decay rate for the stationary distribution of the second queue by controlling the transition structure when the second queue is empty. In Sections 7 and 8, we turn to the question of the decay rate of the probabilities that the process hits level $n$ before level 0. Section 7 deals with general QBD processes; Section 8 deals with the specific case of the tandem Jackson network.

**2. QBD processes.** A level-independent QBD process is a continuous-time Markov chain $(Y_t, J_t, t \geq 0)$ on the state space $\{0, 1, \dots\} \times \{0, \dots, m\}$, whose generator $Q$ has a block tridiagonal representation

$$
(1) \qquad Q = \begin{pmatrix} \widetilde{Q}_1 & Q_0 & & & \\ Q_2 & Q_1 & Q_0 & & \\ & Q_2 & Q_1 & Q_0 & \\ & & Q_2 & Q_1 & Q_0 & \\ & & & \ddots & \ddots & \ddots \end{pmatrix}.
$$

Here, the matrices $Q_0, Q_1, Q_2$ and $\widetilde{Q}_1$ are $(m+1) \times (m+1)$ matrices. The parameter $m$ may be finite or infinite. The random variable $Y_t$ is called the *level* of the process at time $t$ and the random variable $J_t$ is called the *phase* of the process at time $t$.

To avoid complications, we assume that the following condition is satisfied. A discrete-time version of this condition appeared in [9].

CONDITION 2.1. The continuous-time Markov chain on the set $\mathbb{Z} \times \{0, \dots, m\}$, with generator

$$
(2) \qquad \begin{pmatrix} \ddots & \ddots & \ddots & & \\ & Q_2 & Q_1 & Q_0 & \\ & & Q_2 & Q_1 & Q_0 \\ & & & \ddots & \ddots & \ddots \end{pmatrix},
$$

is irreducible.

There are a number of consequences of Condition 2.1, which we shall use later. This condition is satisfied in the QBD process model for the tandem queue which is presented in Section 3.

We set the stage by mentioning some well-known results, at the same time fixing some notation. By Theorem 3.2 of [13], the limiting probabilities $\pi_{kj} := \lim_{t \to \infty} \mathbb{P}(Y_t = k, J_t = j)$ exist. Let us define the vectors $\boldsymbol{\pi}_k = (\pi_{k0}, \ldots, \pi_{km})$, for $k = 0, 1, \ldots,$ and $\boldsymbol{\pi} = (\boldsymbol{\pi}_0, \boldsymbol{\pi}_1, \ldots)$. Then

$$\boldsymbol{\pi}_k = \boldsymbol{\pi}_0 R^k, \qquad k \geq 0, \tag{3}$$

where $R$ is the minimal nonnegative solution to the equation

$$Q_0 + R Q_1 + R^2 Q_2 = 0. \tag{4}$$

The matrix $R$ has a probabilistic interpretation. Let $\mu_i$ be the mean sojourn time in state $(k, i)$, for $k \geq 1$. Then $R(i, j)$ is $\mu_i$ times the total expected time spent in state $(k+1, j)$ before first return to level $k$, starting from state $(k, i)$.

Before turning to the relation between the matrix $R$ and the decay rates of interest, we discuss the issue of ergodicity, both for $m < \infty$ and $m = \infty$, noting a small inaccuracy in the literature regarding the latter case.

THEOREM 2.2.    *The QBD process is ergodic, that is, $\boldsymbol{\pi}$ is positive and has components which sum to unity, if and only if there exists a probability measure $\mathbf{y}_0$ such that*

$$\mathbf{y}_0(\widetilde{Q}_1 + R Q_2) = \mathbf{0} \tag{5}$$

*and*

$$\mathbf{y}_0 \boldsymbol{\nu} < \infty, \tag{6}$$

*where $\boldsymbol{\nu} = (I + R + R^2 + \cdots)\mathbf{1}$. In this case*

$$\boldsymbol{\pi}_0 = \mathbf{y}_0 / \mathbf{y}_0 \boldsymbol{\nu}. \tag{7}$$

The matrix $\widetilde{Q}_1 + R Q_2$ in (5) is the generator of the process of $(Y_t, J_t)$ filtered so that it is observed only when it is in level 0. Thus, the condition that there exists a probability measure satisfying (5) states that the filtered process at level 0 must be ergodic.

In [13], condition (6) is replaced by the elementwise condition

$$\boldsymbol{\nu} < \infty. \tag{8}$$

For the case $m < \infty$, both conditions are equivalent. However, when $m = \infty$, the latter condition is *not* sufficient, since it does not guarantee that $\boldsymbol{\pi}_0$ is nonzero; we may have $\mathbf{y}_0 \boldsymbol{\nu} = \infty$ even when $\boldsymbol{\nu}$ is finite.

Specializing to the case $m < \infty$, inequality (6) is satisfied if and only if

$$(9) \qquad \mathrm{sp}(R) < 1,$$

with $\mathrm{sp}(R)$ denoting the spectral radius of $R$.

If $m$ is finite and there exists a vector $\mathbf{x}$ with $\mathbf{x1} = 1$ such that

$$(10) \qquad \mathbf{x}(Q_0 + Q_1 + Q_2) = 0$$

and

$$(11) \qquad \mathbf{x}Q_0\mathbf{1} < \mathbf{x}Q_2\mathbf{1},$$

then the QBD is positive recurrent; see [11] or [10]. Under the additional assumption that $\widetilde{Q}_1 = Q_1 + Q_2$, this was proved for the infinite case by Tweedie [18]. Equation (11) can be interpreted as requiring that "the average drift of the level process is negative."

We now turn to the decay rate of the stationary distribution, assuming that the QBD process is ergodic. This decay rate is sometimes also referred to as the *caudal characteristic*. We start with a known result for the case $m < \infty$, stating that the geometric decay rate is given by the spectral radius of $R$. In [7], page 205, it was shown that

$$(12) \qquad \lim_{K \to \infty} \frac{\sum_i \pi_{Ki}}{(\mathrm{sp}(R)^K)} = \kappa,$$

where $\kappa$ is a constant. In other words, the marginal stationary probability that the QBD is in level $K$ decays geometrically with rate $\mathrm{sp}(R)$.

Turning to the case $m = \infty$ the situation becomes more complicated, and we come to the core of one of the problems that are dealt with in this paper. We are looking for an "infinite-dimensional" analogue of the limiting result (12), and in particular for the role of the spectral radius of $R$ in it. Clearly, $R$ is now a square matrix of size $\infty$.

There are at least two candidates to consider for this analogue. One approach would be to consider $R$ to be a linear operator from the Banach space $\ell^1$ to itself. We could then hope that the decay rate we are looking for is given by the spectral radius of this operator, if it exists. We shall take a different approach, and use the infinite-dimensional analogue of the Perron–Frobenius eigenvalue of $R$. This is the *convergence norm* of $R$.

Some relevant concepts about the Perron–Frobenius theory of nonnegative matrices are recalled below. For details we refer to [15, 16].

For a finite-dimensional, square, irreducible and nonnegative matrix $A$, there exists a strictly positive eigenvalue which is simple and is greater than or equal to the modulus of all the other eigenvalues. To this eigenvalue corresponds a strictly positive eigenvector. The eigenvalue is called the *Perron–Frobenius eigenvalue* of $A$.

To a large extent, this result can be extended to *infinite-dimensional* matrices. Let $A$ be a nonnegative, aperiodic and irreducible matrix. We would like to prove the existence of a strictly positive $\xi$ and a strictly positive vector $\mathbf{x}$ such that

$$(13) \qquad\qquad \mathbf{x}A = \xi\mathbf{x}.$$

The power series

$$\sum_{k=0}^{\infty} A^k(i,j)z^k$$

has a convergence radius $\alpha$, $0 \leq \alpha < \infty$, *independent* of $i$ and $j$. This common convergence radius is called the *convergence parameter* of the matrix $A$. When $\sum_{k=0}^{\infty} A^k(i,j)\alpha^k$ converges, the matrix is called $\alpha$-*transient*. Otherwise it is called $\alpha$-*recurrent*. An $\alpha$-recurrent matrix $A$ is $\alpha$-*null* if $\lim_{k\to\infty} A^k(i,j)\alpha^k = 0$ and $\alpha$-*positive* otherwise.

The quantity $1/\alpha$ is called the *convergence norm* of $A$. It can be shown to satisfy

$$(14) \qquad\qquad 1/\alpha = \lim_{k\to\infty} (A^k(i,j))^{1/k}$$

independently of $i$ and $j$. This implies, in particular, that if the dimension of $A$ is finite, then the convergence norm is exactly the Perron–Frobenius eigenvalue of $A$ ([16], pages 200 and 201).

For $\beta > 0$, a nonnegative vector $\mathbf{x}$ is called a $\beta$-*subinvariant measure* of $A$ if

$$(15) \qquad\qquad \beta\mathbf{x}A \leq \mathbf{x}$$

and a nonnegative vector $\mathbf{y}$ is called a $\beta$-*subinvariant vector* of $A$ if

$$(16) \qquad\qquad \beta A\mathbf{y} \leq \mathbf{y}.$$

The measure $\mathbf{x}$ and vector $\mathbf{y}$ are called $\beta$-*invariant* when equality holds in (15) and (16), respectively.

The infinite-dimensional analogue of the Perron–Frobenius result is the following (see, e.g., [16], Theorems 6.2 and 6.3):

No $\beta$-subinvariant measure can exist for $\beta > \alpha$. If $A$ is $\alpha$-recurrent, then there exists a strictly positive $\alpha$-invariant measure. If $A$ is $\alpha$-transient, then there exists an $\alpha$-subinvariant measure that is not invariant: there may or may not exist an $\alpha$-invariant measure. By applying the above to the transpose of $A$, similar conclusions can be reached about $\alpha$-invariant vectors.

It is a common misconception to believe that $1/\alpha$ is the largest "eigenvalue" of $A$. This is true in the finite-dimensional case, but not in the infinite-dimensional case. The result above states only that there cannot be any nonnegative $\mathbf{x}$ satisfying (13) with $\xi < 1/\alpha$. In fact, in this paper we shall encounter examples of matrices $A$ such that (13) is satisfied by a positive vector $\mathbf{x}$ for $\xi > 1/\alpha$.

For infinite-dimensional matrices $A$ it is useful to know when the convergence parameter $\alpha$ can be found as a limit of convergence parameters $\{\alpha^{(k)}\}$ from a sequence $\{A^{(k)}\}$ of finite-dimensional matrices. For example, in Theorem 6.8 of [16], it is shown that the convergence parameters of the $(n \times n)$ northwest corner truncations of $A$ converge to the convergence parameter of $A$. The following result will be of use to us in Sections 7 and 8.

LEMMA 2.3. *Let $\{A^{(k)}\}$ be a sequence of nonnegative matrices that increases elementwise to an irreducible matrix $A$, as $k \to \infty$. Let $\alpha^{(k)}$ denote the convergence parameter of $A^{(k)}$ and let $\alpha$ be the convergence parameter of $A$. Then the sequence $\alpha^{(k)}$ is decreasing with $\lim_{k\to\infty} \alpha^{(k)} = \alpha$.*

PROOF. By (14) and the fact that $A^{(k)}(i,j) \le A^{(k+1)}(i,j)$, we have

$$\alpha^{(k+1)} \le \alpha^{(k)}.$$

Therefore the sequence $\{\alpha^{(k)}\}$ is decreasing and its limit $\alpha^{(\infty)}$ must exist. Also, by (14) and the fact that $A^{(k)}(i,j) \le A(i,j)$ for all $k$, we have

$$\alpha \le \alpha^{(\infty)}.$$

Now let $\mathbf{y}^{(k)} = (y_1^{(k)}, y_2^{(k)}, \dots)^T$ be an $\alpha^{(k)}$-subinvariant vector of $A^{(k)}$, with $y_1^{(k)} = 1$, and let $\mathbf{y}^* = \liminf_k \mathbf{y}^{(k)}$, elementwise. Then we know that

$$\alpha^{(k)} A^{(k)} \mathbf{y}^{(k)} \le \mathbf{y}^{(k)}.$$

Taking $\liminf_{k\to\infty}$ of both sides and using Fatou's lemma, we have, for each $i$,

$$(17) \qquad \alpha^{(\infty)} \sum_{j=1}^{\infty} A(i,j) y_j^* \le y_i^*.$$

Iterating this, we find that, for $\nu \ge 1$,

$$\alpha^{(\infty)\nu} \sum_{j=1}^{\infty} A^\nu(i,j) y_j^* \le y_i^*.$$

Since $y_1^* = 1$ and $A$ is irreducible, this shows that $y_j^* < \infty$ for all $j$ and, by (17), that $\mathbf{y}^*$ is an $\alpha^{(\infty)}$-subinvariant vector of $A$. Since no $\beta$-subinvariant vector can exist for $\beta > \alpha$, we must have $\alpha^{(\infty)} \le \alpha$ and thus $\alpha^{(\infty)} = \alpha$. □

Now let us turn back to the problem of determining the decay rate of $R^n$. From the definition, it follows that if $\sum_{n=0}^{\infty} R^n(i,j)$ is convergent for all $i$ and $j$, then the convergence norm of $R$ must be less than or equal to 1. It is thus tempting to think that the decay rate of the stationary distribution must be given by the convergence norm. However, we have to be careful. As we noted above, it is a common misconception to believe that the convergence norm is the largest "eigenvalue."

Assume that $\mathbf{w}$ is a $z^{-1}$-invariant measure of $R$ such that $\mathbf{w}\sum_{i=0}^{\infty} R^i$ is finite. Then $z$ must be less than 1. To see this, note that the monotone convergence theorem implies that $\mathbf{w}\sum_{i=0}^{k} R^i$ converges elementwise to $\mathbf{w}\sum_{i=0}^{\infty} R^i$. This means that $\mathbf{w}\sum_{i=0}^{k} z^i$ converges elementwise to a finite vector, which can be the case only when $z < 1$. This leads to the following result.

THEOREM 2.4. *Consider an irreducible QBD process with a finite or infinite phase space. If there exists a nonnegative vector $\mathbf{w} \in \ell^1$ and a nonnegative number $z < 1$ such that*

$$(18) \qquad\qquad \mathbf{w}(\widetilde{Q}_1 + RQ_2) = \mathbf{0}$$

*and*

$$(19) \qquad\qquad \mathbf{w}R = z\mathbf{w},$$

*then the QBD is ergodic, and, for all fixed $i = 0, 1, \ldots,$*

$$(20) \qquad\qquad \frac{\pi_{Ki}}{z^K} = \mathbf{w}_i$$

*for all $K$.*

Theorem 2.4 shows that if $\boldsymbol{\pi}_0$ is a $z^{-1}$-invariant measure of $R$ for some $z$, then the stationary distribution of $(Y_t, J_t)$ has the *level-phase independence property* (see [8]) and decays at rate $z$. If $\boldsymbol{\pi}_0$ is a (finite) linear combination of more than one $\mathbf{w}$ (not necessarily nonnegative) that satisfies $\mathbf{w}R = \xi\mathbf{w}$ for some $\xi$, then the stationary distribution does not have this property. The decay rate is then given by the value of $\xi$ in the linear combination for which $|\xi|$ is the largest.

For the case $m < \infty$, any $\boldsymbol{\pi}_0$ is a finite linear combination of eigenvectors of $R$. The left eigenvector of the eigenvalue $\mathrm{sp}(R)$ must always be in this linear combination. To see this, recall that $\boldsymbol{\pi}_0$ must be positive, and the Perron–Frobenius right-eigenvector $\mathbf{v}$ of $R$ must be nonnegative and nonzero, which implies that $\boldsymbol{\pi}_0\mathbf{v} > 0$. Now write

$$\boldsymbol{\pi}_0 = \sum_{i=1}^{m} a_i\mathbf{w}_i,$$

where $\mathbf{w}_1$ is the Perron–Frobenius eigenvector. Since, for $i \geq 2$, the $\mathbf{w}_i$ correspond to eigenvalues of $R$ distinct from $\mathrm{sp}(R)$, we know that $\mathbf{w}_i \mathbf{v} = 0$ for these values of $i$. Therefore $\boldsymbol{\pi}_0 \mathbf{v} = a_1 \mathbf{w}_1 \mathbf{v}$, which shows that $a_1 \neq 0$. This explains why the stationary distribution decays at rate $\mathrm{sp}(R)$ when $m$ is finite.

As a final topic in this section, we quote a result that helps us to determine $z^{-1}$-invariant measures of $R$. In the general case $m \leq \infty$ it is easy to see from (4) that, if the row vector $\mathbf{w}$ and scalar $z$ satisfy $\mathbf{w}R = z\mathbf{w}$, then

$$(21) \qquad \mathbf{w}(Q_0 + zQ_1 + z^2Q_2) = \mathbf{0}$$

whenever the change of order of summation involved in using the associative law of matrix multiplication is permitted.

More important, under certain conditions the converse is true as well, again irrespective of whether $m < \infty$ or $m = \infty$. This is shown in the next theorem, which is a statement of Theorem 5.4 of Ramaswami and Taylor [14]. See also [4] for a more detailed analysis for the case where $m < \infty$.

THEOREM 2.5. *Consider a continuous-time ergodic QBD process with generator of the form* (1). *Let* $q_k = -Q_1(k,k)$. *If the complex variable* $z$ *and the vector* $\mathbf{w} = \{w_k\}$ *are such that* $|z| < 1$ *and* $\sum_k |w_k| q_k < \infty$, *then* (21) *implies* (19).

Specializing to the finite case $m < \infty$, we obtain the following corollary.

COROLLARY 2.6. *For an irreducible QBD process with a finite phase space, satisfying* (11), *the eigenvalues of* $R$ *are all the zeros of the polynomial*

$$(22) \qquad \det(Q_0 + zQ_1 + z^2Q_2)$$

*that lie strictly within the unit circle.*

PROOF. From the discussion above it follows that all eigenvalues of $R$ lie within the unit circle. Each such eigenvalue $z$ with corresponding left eigenvector $\mathbf{w}$ satisfies (21). Conversely, by Theorem 2.5, all solutions $(z, \mathbf{w})$ to (21) with $z$ within the unit circle must be eigenvalue–left-eigenvector pairs of $R$, because the condition $\sum_k |w_k| q_k < \infty$ is automatically satisfied. In particular, the eigenvalues of $R$ are the zeros of (22) within the unit circle. □

**3. The tandem Jackson network seen as a QBD process.** We now turn to a specific class of QBD processes which may have infinitely many phases. It models a simple Jackson network consisting of two queues in tandem (see Figure 1). Customers arrive at the first queue according to a Poisson process with rate $\lambda$. The service time of customers at the first queue is exponentially distributed with parameter $\mu_1$. On leaving the first queue, customers enter the second queue, where their service time has an exponential distribution with parameter $\mu_2$. The capacity of the first queue is denoted by $m$, which may be finite or infinite. In the case when $m$ is finite, customers that arrive to find the first queue full are rejected. For $i = 1, 2$ let

$$\rho_i = \frac{\lambda}{\mu_i}$$

and let $J_t$ and $Y_t$ denote the number of customers in the first and second queue at time $t$, respectively.

We shall examine the behavior of the two-dimensional Markov chain $(Y_t, J_t)$, viewed as a QBD process in which $Y_t$ represents the level and $J_t$ represents the phase. The transition intensities of this QBD process are depicted in Figure 2. When the capacity of the first queue is infinite, the phase space of this QBD process is infinite and the boundary denoted by $m$ in the figure is not present.

For the case where $m < \infty$, the $(m + 1) \times (m + 1)$-matrices $Q_0, Q_1, Q_2$ and $\widetilde{Q}_1$ in (1) are given by

$$Q_0 = \begin{pmatrix} 0 & \dots & & & \\ \mu_1 & 0 & \dots & & \\ & \mu_1 & 0 & \dots & \\ & & \ddots & \ddots & \\ & & & \mu_1 & 0 \end{pmatrix}, \qquad Q_2 = \begin{pmatrix} \mu_2 & & & & \\ & \mu_2 & & & \\ & & \mu_2 & & \\ & & & \ddots & \\ & & & & \mu_2 \end{pmatrix},$$

$$Q_1 = \begin{pmatrix} -(\lambda+\mu_2) & \lambda & & & \\ & -(\lambda+\mu_1+\mu_2) & \lambda & & \\ & & -(\lambda+\mu_1+\mu_2) & \lambda & \\ & & & \ddots & \ddots \\ & & & & -(\mu_1+\mu_2) \end{pmatrix}$$

and

$$\widetilde{Q}_1 = \begin{pmatrix} -\lambda & \lambda & & & \\ & -(\lambda+\mu_1) & \lambda & & \\ & & -(\lambda+\mu_1) & \lambda & \\ & & & \ddots & \ddots \\ & & & & -\mu_1 \end{pmatrix}.$$

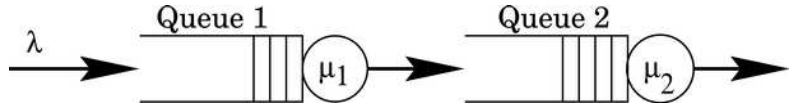

Fig. 1.   *A tandem Jackson network.*

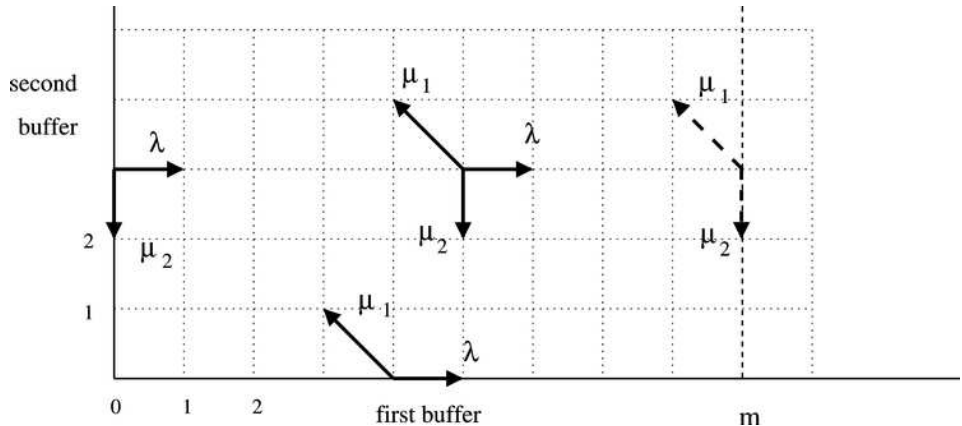

Fig. 2.   *The transition intensities for the tandem network.*

Obviously, Condition 2.1 is satisfied in this case and the stability condition (11) translates into

$$(23) \qquad \rho_2 < \frac{1-\rho_1^{m+1}}{1-\rho_1^m}, \qquad \rho_1 \neq 1,$$

$$(24) \qquad \rho_2 < 1 + \frac{1}{m}, \qquad \rho_1 = 1.$$

For the case where $m = \infty$, the tridiagonal blocks are given by the infinite-dimensional matrices

$$Q_0 = \begin{pmatrix} 0 & \dots & & \\ \mu_1 & 0 & \dots & \\ & \mu_1 & 0 & \dots \\ & & \ddots & \ddots \end{pmatrix}, \qquad Q_2 = \begin{pmatrix} \mu_2 & & & \\ & \mu_2 & & \\ & & \mu_2 & \\ & & & \ddots \end{pmatrix},$$

$$Q_1 = \begin{pmatrix} -(\lambda+\mu_2) & \lambda & & \\ & -(\lambda+\mu_1+\mu_2) & \lambda & \\ & & -(\lambda+\mu_1+\mu_2) & \lambda \\ & & & \ddots & \ddots \end{pmatrix}$$

and

$$\widetilde{Q}_1 = \begin{pmatrix} -\lambda & \lambda & & \\ & -(\lambda + \mu_1) & \lambda & \\ & & -(\lambda + \mu_1) & \lambda & \\ & & & \ddots & \ddots \end{pmatrix}.$$

In this case, the well-known condition under which both queues are stable is

$$(25) \qquad\qquad\qquad \lambda < \min\{\mu_1, \mu_2\}.$$

For both finite and infinite $m$, we are interested in the decay rate of the stationary distribution of the tandem network as the number in the second queue becomes large and its relation to the spectral properties of the matrix $R$. It will be convenient to index this matrix with the size of the waiting room at the first queue. Thus we shall write $R_m$ for the situation where the size of this waiting room is $m$ and, in particular, $R_\infty$ when the waiting room at the first queue is unlimited.

For infinite $m$, under condition (25), it follows from the results of Burke [1] that the arrival process to the second queue is a Poisson process with parameter $\lambda$ and so the second queue behaves like an $M/M/1$ queue with arrival rate $\lambda$ and service rate $\mu_2$. Thus the stationary distribution of the second queue is geometric with parameter $\rho_2$ and its decay rate is simply $\rho_2$. However, it is not at all clear how this decay rate corresponds to the spectral properties of the infinite-dimensional matrix $R_\infty$.

To study the spectral properties of $R_\infty$ we shall make use of Theorem 2.5. To facilitate our development, we introduce some notation.

For each $z$ with $|z| < 1$, $z \neq 0$, let $Q(z)$ be the infinite-dimensional tridiagonal matrix $(Q_0 + zQ_1 + z^2Q_2)/z$, that is,

$$(26) \quad Q(z) = \begin{pmatrix} -\lambda - \mu_2 + \mu_2 z & \lambda & & \\ \mu_1/z & -\lambda - \mu_1 - \mu_2 + \mu_2 z & \lambda & \\ & \ddots & \ddots & \ddots \end{pmatrix},$$

and let $Q^{(n)}(z)$ denote the $(n \times n)$ northwest corner truncation of $Q(z)$.

For the finite case, define the $(n \times n)$-matrix $\widehat{Q}^{(n)}(z)$ as

$$(27)$$
$$\widehat{Q}^{(n)}(z)$$
$$= \begin{pmatrix} -\lambda - \mu_2 + \mu_2 z & \lambda & & \\ \mu_1/z & -\lambda - \mu_1 - \mu_2 + \mu_2 z & \lambda & \\ & \ddots & \ddots & \ddots \\ & & \mu_1/z & -\mu_1 - \mu_2 + \mu_2 z \end{pmatrix}.$$

The significance of these matrices follows from Theorem 2.5. The infinite-dimensional row vector $\mathbf{w}$ satisfies $\mathbf{w}R_\infty = z\mathbf{w}$ for $z \neq 0$ with $|z| < 1$, if $\sum_k |w_k|q_k < \infty$ and $\mathbf{w}$ satisfies

$$\mathbf{w}Q(z) = \mathbf{0}. \tag{28}$$

For the tandem queue, $q_k$ is constant for $k \geq 1$ and so the condition that $\sum_k |w_k|q_k < \infty$ is equivalent to requiring that $\mathbf{w} \in \ell^1$.

For the case $m < \infty$, the $(m+1)$-dimensional row vector $\mathbf{w}$ is a left eigenvector of $R_m$ corresponding to eigenvalue $z \neq 0$ with $|z| < 1$, if and only if it satisfies

$$\mathbf{w}\widehat{Q}^{(m+1)}(z) = \mathbf{0}. \tag{29}$$

REMARK 3.1. Readers may note that (28) and (29) are not exactly equivalent to (21). The latter follow from the former only if $z \neq 0$. In fact, for the tandem network model, the vector $(1, 0, 0, \dots)$ satisfies (21) with $z = 0$.

By using the physical interpretation of $R_m$, we can see that the interesting $z^{-1}$-invariant measures of $R_m$ are the ones for which (28) and (29) are satisfied. For the tandem Jackson network, the expected time spent in any state at level $k+1$ before the process returns to level $k$ is nonzero if the process starts in a state $(k, i)$ with $i > 0$. Thus we know immediately from its physical interpretation that $R_m(i, j)$ is strictly positive for all $i \geq 1$ and $j \geq 0$. On the other hand, it is impossible to visit level $k+1$ starting in state $(k, 0)$ without visiting a state $(k, i)$ with $i \geq 1$ first, and so $R_m(0, j) = 0$ for all $j$. Thus $R_m$ decomposes its indices into two communicating classes, $\mathcal{C}_1 \equiv \{0\}$ and $\mathcal{C}_2 \equiv \{1, 2, \dots\}$. The eigenvector $(1, 0, 0, \dots)$ of $R_m$ with corresponding eigenvalue 0 has support on $\mathcal{C}_1$. All other $z^{-1}$-invariant measures of $R_m$ have the form $(w_0, \mathbf{w}_1)$, where $\mathbf{w}_1$ is a $z^{-1}$-invariant measure of the positive submatrix $\widetilde{R}_m$ corresponding to $\mathcal{C}_2$. These are the $\mathbf{w}$ and $z$ for which (28) and (29) are satisfied.

**4. The case where $m$ is infinite.** Before we start studying (28), we first give some preliminaries. In this and the following sections, we shall frequently use the function

$$\tau(z) \equiv -\lambda - \mu_1 - \mu_2(1 - z) + 2\sqrt{\frac{\lambda\mu_1}{z}}. \tag{30}$$

It is easy to see that $\tau(z)$ is convex on $(0, 1)$ with $\lim_{z \to 0} \tau(z) = \infty$, and $\tau(1) = -(\sqrt{\lambda} - \sqrt{\mu_1})^2$. Thus there is a unique value $\eta \in (0, 1)$ with $\tau(\eta) = 0$, and, for $z \in (0, 1)$, $\tau(z) < 0$ if and only if $z > \eta$.

We shall also frequently refer to the relationship between $\eta$, $\rho_1$ and $\rho_2$ in the respective cases when $\mu_1 \leq \mu_2$ and $\mu_1 > \mu_2$. These are summarized in the following lemma.

LEMMA 4.1.   (a) *When $\mu_1 \leq \mu_2$, $0 < \eta \leq \rho_2 \leq \rho_1 < 1$.*
(b) *When $\mu_1 > \mu_2$, $0 < \rho_1 < \eta < \rho_2 < 1$.*

PROOF.   Observe that $\tau(\rho_2) \leq 0$, which immediately gives us that $\rho_2 \geq \eta$, and $\tau(\rho_1) = (1 - \lambda/\mu_1)(\mu_1 - \mu_2)$, which gives us that $\rho_1 \geq \eta$ when $\mu_1 \leq \mu_2$ and $\rho_1 < \eta$ when $\mu_1 > \mu_2$. Together with the fact that $\rho_1 \geq \rho_2$ if and only if $\mu_1 \leq \mu_2$, this proves the lemma.   □

Now consider the system of equations (28) where $z$ is fixed such that $z \in (-1, 1)$, $z \neq 0$. Writing out the system, we have

$$(31) \qquad -(\lambda + \mu_2(1-z))zw_0 + \mu_1 w_1 = 0,$$

$$(32) \qquad \lambda z w_{k-1} - (\lambda + \mu_1 + \mu_2(1-z))zw_k + \mu_1 w_{k+1} = 0, \qquad k \geq 1.$$

After substituting $w_k = u^k$ in (32), we derive the characteristic equation,

$$(33) \qquad \mu_1 u^2 - (\lambda + \mu_1 + \mu_2(1-z))zu + \lambda z = 0.$$

Since the discriminant of (33) is positive if and only if $z < 0$ or $\tau(z) < 0$, the form of the solution now depends on the location of $z$ relative to $0$ and $\eta$. We proceed by giving the solution for $w_k$ in the cases $-1 < z < 0$ and $\eta < z < 1$. This is

$$(34) \qquad w_k = c_1 u_1^k + c_2 u_2^k,$$

where

$$(35) \qquad u_{1,2} = \frac{(\lambda + \mu_1 + \mu_2(1-z))z \pm \sqrt{(\lambda + \mu_1 + \mu_2(1-z))^2 z^2 - 4\lambda \mu_1 z}}{2\mu_1}.$$

The coefficients $c_1$ and $c_2$ can be derived from

$$(36) \qquad c_1 + c_2 = 1,$$

$$(37) \qquad c_1 u_1 + c_2 u_2 = \frac{1}{\mu_1}(\lambda + \mu_2(1-z))z,$$

where the first equation is due to the (arbitrary) normalizing assumption that $w_0 = 1$, and the second equation follows from boundary equation (31). Thus, we find

$$c_{1,2} = \frac{1}{2} \pm \frac{(\lambda - \mu_1 + \mu_2(1-z))z}{2\sqrt{(\lambda + \mu_1 + \mu_2(1-z))^2 z^2 - 4\lambda \mu_1 z}}.$$

When $z = \eta$, the vector $\mathbf{w}$ is given by

$$(38) \qquad w_k = u^k(1 + ck),$$

with $u = \sqrt{\rho_1 \eta}$ and $c = 1 - \sqrt{\eta/\rho_1}$, while for $0 < z < \eta$ the real solution is given by

$$(39) \qquad w_k = (\cos(k\phi) + c\sin(k\phi))|u|^k,$$

with $|u| = \sqrt{\rho_1 z}$,

$$\phi = \arctan\left(\frac{\sqrt{4\lambda\mu_1 z - (\lambda + \mu_1 + \mu_2(1-z))^2 z^2}}{(\lambda + \mu_1 + \mu_2(1-z))z}\right)$$

and

$$c = \frac{(\lambda + \mu_2(1-z))\sqrt{z/\lambda\mu_1} - \cos(\phi)}{\sin(\phi)}.$$

As we pointed out after equation (28), in order to use Theorem 2.5 to establish whether $\mathbf{w}$ is indeed a $z^{-1}$-invariant measure of $R$, we need to verify whether $\mathbf{w} \in \ell^1$.

LEMMA 4.2. *The vector* $\mathbf{w}$ *is an element of* $\ell^1$ *if and only if*

$$z_1 < z < \mu_1/\mu_2,$$

*where* $z_1 = (2\lambda + \mu_1 + \mu_2 - \sqrt{(2\lambda + \mu_1 + \mu_2)^2 + 4\mu_1\mu_2})/(2\mu_2) < 0.$

PROOF. First note that, for $0 < z \le \eta$, the form of (38) and (39) shows that it is certain that $\mathbf{w} \in \ell^1$. Thus we need only consider the case when the roots $u_1$ and $u_2$ are real. This occurs when $-1 < z < 0$ or $\eta < z < 1$.

Unless $z = 1$ or $z = \rho_2$ both $c_1$ and $c_2$ are nonzero, so for $\mathbf{w}$ to be in $\ell^1$ it is necessary and sufficient that both $u_1$ and $u_2$ are in $(-1, 1)$. To study when this is the case, let $f(u)$ be the left-hand side of (33). Then the statement that the roots $u_1$ and $u_2$ are in $(-1, 1)$ is equivalent to saying that both $f(-1) > 0$, $f(1) > 0$, $f'(-1) < 0$ and $f'(1) > 0$.

When $-1 < z < 0$, $f(1)$ is always positive and the condition that $f(-1)$ is positive reduces to $\mu_1 + (2\lambda + \mu_1 + \mu_2(1-z))z > 0$, which is the same as saying that $z_1 < z$. Furthermore, $f'(1)$ is always positive and $f'(-1)$ can be written as $-f(1) - \mu_1 + \lambda z$, which is negative if $f(1)$ is positive. When $\eta < z < 1$, $f(-1)$ is always positive and the condition that $f(1)$ is positive is $(1 - z)(\mu_1 - z\mu_2) > 0$, which reduces to $z < \mu_1/\mu_2$. Furthermore $f'(1)$ is always negative and $f'(1) = f(1) - \lambda + \mu_1$, which is positive when $f(1)$ is positive by the stability condition (25).

The observations that $\eta < \mu_1/\mu_2$, established by verifying that $\tau(\mu_1/\mu_2) < 0$, and $-1 < z_1 < 0$ complete the proof. $\square$

COROLLARY 4.3. *When $\mu_1 \leq \mu_2$, the system of equations*

$$(40) \qquad \mathbf{w}R_\infty = z\mathbf{w}$$

*has solutions $\mathbf{w} \in \ell^1$ for all $z \in (z_1, \mu_1/\mu_2)$.*

*When $\mu_1 > \mu_2$, the system (40) has solutions $\mathbf{w} \in \ell^1$ for all $z \in (z_1, 1)$.*

Note that if $\mu_1 < \mu_2$, it is not certain whether (40) has solutions for $z \in [\mu_1/\mu_2, 1)$, but any such solutions will not be in $\ell^1$. In Remark 6.3 we show that such solutions exist only for $z = \mu_1/\mu_2$.

For us to be able to apply Theorem 2.4, the vector $\mathbf{w}$ must be nonnegative. To investigate this, we start by generalizing (31) and (32) to

$$(41) \qquad P_0(x; z) = 1,$$

$$(42) \qquad \frac{\mu_1}{z} P_1(x; z) = x + \lambda + \mu_2(1 - z),$$

$$\frac{\mu_1}{z} P_n(x; z) = (x + \lambda + \mu_1 + \mu_2(1 - z))P_{n-1}(x; z) - \lambda P_{n-2}(x; z),$$
$$(43)$$
$$n \geq 2.$$

For any given real and positive value of $z$, (41)–(43) define a sequence of *orthogonal polynomials* $P_n(x; z)$. When $x = 0$, they reduce to (31) and (32), from which we deduce the fact that $w_n = P_n(0; z)$. Moreover, we shall see that $P_n(0; z)$ is positive for all $n$ if and only if the zeros of all the $P_n(x; z)$ are less than zero. Thus we can study conditions for the positivity of $\mathbf{w}$ via the properties of the polynomials $P_n(x; z)$.

LEMMA 4.4. *For $z > 0$, the sequence $\{P_n(x; z)\}$ satisfies the orthogonality relationship*

$$\int_{\mathrm{supp}(\psi)} P_n(x; z) P_m(x; z) \psi(dx) = \left(\frac{z\lambda}{\mu_1}\right)^n \delta_{n,m},$$

*where*

$$\mathrm{supp}(\psi) = \begin{cases} [\sigma(z), \tau(z)], & \text{if } z \leq \rho_1, \\ [\sigma(z), \tau(z)] \cup \{\chi(z)\}, & \text{if } z > \rho_1, \end{cases}$$

$\tau(z)$ *is given by (30),*

$$(44) \qquad \sigma(z) = -\lambda - \mu_1 - \mu_2(1 - z) - 2\sqrt{\frac{\lambda\mu_1}{z}}$$

*and*

$$(45) \qquad \chi(z) = \left(\frac{\lambda}{z} - \mu_2\right)(1 - z).$$

*The measure $\psi$ is given by*

$$\psi(dx) = \frac{2}{\pi} \frac{\sqrt{1 - (x + \lambda + \mu_1 + \mu_2(1-z))^2 z / 4\lambda\mu_1}}{1 - (x + \lambda + \mu_2(1-z))z/\lambda} \, dx, \qquad \sigma \le x \le \tau,$$

$$\psi(\{\chi(z)\}) = 1 - \frac{\lambda}{z\mu_1} \qquad \textit{if } z > \rho_1.$$

PROOF. For fixed $z > 0$, let

$$(46) \qquad T_n(x) = \left(\sqrt{\frac{\mu_1}{z\lambda}}\right)^n P_n\left(2x\sqrt{\frac{\lambda\mu_1}{z}} - \lambda - \mu_1 - \mu_2(1-z); z\right).$$

It follows that $T_0(x) = 1$, $T_1(x) = 2x - b$ and $T_n(x) = 2xT_{n-1}(x) - T_{n-2}(x)$, where

$$b = \sqrt{\frac{z\mu_1}{\lambda}}.$$

The $T_n$'s are *perturbed Chebyshev polynomials*, for which the orthogonalizing relationship is given (see [2], pages 204 and 205) by

$$\frac{2}{\pi} \int_{-1}^{1} T_n(x) T_m(x) \frac{\sqrt{1 - x^2}}{1 + b^2 - 2bx} \, dx$$

$$+ \mathbf{1}_{\{|b|>1\}} T_n\left(\frac{b}{2} + \frac{1}{2b}\right) T_m\left(\frac{b}{2} + \frac{1}{2b}\right)\left(1 - \frac{1}{b^2}\right) = \delta_{n,m},$$

where $\mathbf{1}_{\{|b|>1\}} = 1$ if $|b| > 1$ and 0 otherwise. Substituting (46) and rewriting yields the result. $\square$

As a consequence we have the following.

LEMMA 4.5. *For each value of $z > 0$, $P_n(x; z)$ has $n$ distinct real zeros $x_{n,1} < \cdots < x_{n,n}$ and these zeros interlace. That is, for all $n \ge 2$ and $i = 1, \ldots, n-1$,*

$$x_{n,i} < x_{n-1,i} < x_{n,i+1}.$$

PROOF. The lemma follows from a well-known result for orthogonal polynomial sequences (see [2], Theorem 5.3). $\square$

The support of the measure $\psi$ is intimately related to the limiting behavior of the zeros of the $P_n(x; z)$. Some results are stated in the lemma below.

Lemma 4.6.   *The sequences of smallest, second-largest and largest zeros of the $P_n(x; z)$ possess the following properties:*

   $\{x_{n,1}\}_{n=1}^{\infty}$ *is a strictly decreasing sequence with limit* $\sigma(z)$;

   $\{x_{n,n-1}\}_{n=1}^{\infty}$ *is a strictly increasing sequence with limit* $\tau(z)$;

   $\{x_{n,n}\}_{n=1}^{\infty}$ *is a strictly increasing sequence with limit* $\chi_1(z)$,

*where*

$$\chi_1(z) = \sup(\text{supp}(\psi)) = \begin{cases} \tau(z), & \text{if } z \leq \rho_1, \\ \chi(z), & \text{if } z > \rho_1. \end{cases}$$

For a proof, see [2], Section II.4.

Lemma 4.7.   *Let $z > 0$. Then $P_n(x; z)$ is positive for all $n$ if and only if $x \geq \chi_1(z)$.*

Proof.   The leading coefficient of $P_n(x; z)$ is positive for all $n$, which implies that $P_n(x; z)$ is positive for $x > x_{n,n}$. Since $x_{n,n}$ is strictly increasing, we know that $P_n(x; z)$ is positive for all $n$ if $x \geq \chi_1(z)$. Conversely, $P_k(x; z)$ is negative for $x \in (x_{k-1,k}, x_{k,k})$ and so the interleaving property given in Lemma 4.5 implies that, for every $x < x_{n,n}$, $P_k(x; z)$ is less than zero for at least one $k \in \{1, \dots, n\}$. Thus, if $x < \chi_1(z)$, $P_k(x; z)$ is less than zero for at least one $k \in Z_+$.   □

Next, let us return to the question of when the vector **w** which solves (31) and (32) is positive.

Lemma 4.8.   *The vector **w** is positive if and only if $\chi_1(z) \leq 0$.*

Proof.   This follows immediately from Lemma 4.7 and the fact that, for a given value of $z$, $w_n = P_n(0; z)$.   □

Lemma 4.8 implies that, to decide whether **w** is positive, it is important to know for which values of $z$ the corresponding $\chi_1(z)$ is less than or equal to 0. Since

$$\chi_1(z) = \max(\chi(z), \tau(z)) = \begin{cases} \tau(z), & \text{for } z \leq \rho_1, \\ \chi(z), & \text{for } z > \rho_1, \end{cases}$$

the statement that $\chi_1(z) \leq 0$ implies that, for $z \leq \rho_1$,

$$\tau(z) \leq 0 \quad \text{and so} \quad z \geq \eta$$

and, for $z > \rho_1$,

$$\chi(z) \leq 0 \quad \text{and so} \quad z \geq \rho_2.$$

When $\mu_1 \leq \mu_2$, we know from Lemma 4.1 that

$$\eta \leq \rho_2 \leq \rho_1$$

and so $\mathbf{w}$ is positive for all $z \in [\eta, 1)$. When $\mu_1 > \mu_2$, Lemma 4.1 tells us that

$$\rho_1 < \eta < \rho_2.$$

Thus $\mathbf{w}$ is positive only for $z \in [\rho_2, 1)$.

Summarizing this and Corollary 4.3, we have the following theorem.

THEOREM 4.9. *When $\mu_1 \leq \mu_2$, the system of equations* (19) *has positive solutions $\mathbf{w} \in \ell^1$ for all $z \in [\eta, \mu_1/\mu_2)$.*

*When $\mu_1 > \mu_2$, the system* (19) *has positive solutions $\mathbf{w} \in \ell^1$ for all $z \in [\rho_2, 1)$.*

Theorem 4.9 states a very interesting result. Together with Theorem 2.4, it indicates that it might be possible to have level-phase independent stationary distributions of the tandem queue for a range of different $z$. The key point is whether the vector $\mathbf{w}$ that satisfies (19) also satisfies (18).

In fact it has been well-known since the work of Burke [1] and Jackson [5] that the decay rate of the stationary number of customers in the second queue is $\rho_2$ irrespective of whether $\mu_1 \leq \mu_2$ or $\mu_1 > \mu_2$, and not any of the other possible values of $z$. Why should this be the case? The answer is that $\boldsymbol{\pi}_0$, the distribution of $J$ at level 0 satisfying (18), is precisely the vector $\mathbf{w}$ that satisfies (19) with $z = \rho_2$. In other words, the decay rate is *that* value $z$ for which $R_\infty$ which has the proper $z^{-1}$-invariant measure.

This leads us to ask the question that if we varied $\widetilde{Q}_1$, and thus (18), can we get a vector $\mathbf{w}$ that satisfies (19) for a value of $z \neq \rho_2$. If we can do this, we shall have changed the decay rate of the stationary distribution of the number in the second queue by changing the transition structure only when the second queue is empty. In Section 6, we shall see that it is indeed possible to do this.

Before we move on, we shall briefly discuss how the results of Takahashi, Fujimoto and Makimoto [17] apply to the tandem network example. An application of Corollary 1 of [17] shows that if there exists a scalar $z$ and vector $\mathbf{w} \in \ell^1$ that satisfy (28) and a vector $\mathbf{y}$ that satisfies

$$(47) \qquad\qquad Q(z)\mathbf{y} = \mathbf{0}$$

with $\mathbf{w}\mathbf{y} < \infty$ and $z^{-1}\mathbf{w}A_0\mathbf{y} \neq z\mathbf{w}A_2\mathbf{y}$, then (19) is satisfied, $R$ is $z$-positive and the right eigenvector of $R$ is dominated elementwise by $\mathbf{y}$. Furthermore, by Corollary 2 of [17], if $\boldsymbol{\pi}_0$ is such that $\boldsymbol{\pi}_0\mathbf{y} < \infty$, then the decay rate of the QBD process is equal to $z$.

After some calculation, we see that, when $\mu_1 > \mu_2$, the conditions of Corollary 1 of [17] are satisfied with $z = \rho_2$, $\mathbf{w}$ such that $w_k = \rho_1^k$ and $\mathbf{y}$ such that $y_k = \rho_2^{-k}$. Because $\boldsymbol{\pi}_0 = \mathbf{w}$, we can then derive the fact that the decay rate is $\rho_2$. However, this reasoning does not work if $\mu_1 \leq \mu_2$ and, even if $\mu_1 > \mu_2$, by altering $\tilde{Q}_1$, we can create the situation where $\boldsymbol{\pi}_0$ is a $z^{-1}$-invariant measure of $R$ for a different value of $z$. In this case $\boldsymbol{\pi}_0$ must necessarily be such that $\boldsymbol{\pi}_0 \mathbf{y} = \infty$. We give an example of such a construction in Section 6.

**5. The case where $m$ is finite.** In the case where $m$ is finite, because the tandem queue is assumed to be stable, we know by Corollary 2.6 that the nonzero eigenvalues of $R_m$ are given by the values of $z$ within the unit circle for which $\det \hat{Q}^{(m+1)}(z) = 0$. Thus $(z, \mathbf{w})$ is an eigenvalue–eigenvector pair of $R_m$ if and only if zero is an eigenvalue of $\hat{Q}^{(m+1)}(z)$ with corresponding eigenvector $\mathbf{w}$. In the first part of this section, we shall explore the relationship between the values of $x$ for which $\det(xI_{m+1} - \hat{Q}^{(m+1)}(z)) = 0$ and the zeros of a sequence of orthogonal polynomials closely related to the $P_n(x; z)$.

Let the sequence of polynomials $\hat{P}_n(x; z)$, be defined such that $\hat{P}_0(x; z) = 1$ and, for $n \geq 1$,

$$\hat{P}_n(x; z) = P_n(x; z) - \frac{\lambda z}{\mu_1} P_{n-1}(x; z).$$

The polynomials $\hat{P}_n(x; z)$ satisfy the recursion

$$(48) \qquad \hat{P}_0(x; z) = 1,$$

$$(49) \qquad \frac{\mu_1}{z} \hat{P}_1(x; z) = x + \mu_2(1 - z),$$

$$(50) \qquad \frac{\mu_1}{z} \hat{P}_2(x; z) = (x + \lambda + \mu_1 + \mu_2(1 - z))\hat{P}_1(x; z) - \lambda(1 - z),$$

$$\frac{\mu_1}{z} \hat{P}_n(x; z) = (x + \lambda + \mu_1 + \mu_2(1 - z))\hat{P}_{n-1}(x; z) - \lambda \hat{P}_{n-2}(x; z),$$
$$(51)$$
$$n \geq 3.$$

LEMMA 5.1. *For each value of $z > 0$, $\hat{P}_n(x; z)$ has $n$ distinct real zeros $\hat{x}_{n,1} < \cdots < \hat{x}_{n,n}$ which interlace. Moreover, $\hat{x}_{n,n} > x_{n,n}$ and*

$$(52) \qquad x_{n,i} < \hat{x}_{n,i} < x_{n,i+1}, \qquad i = 1, \ldots, n-1.$$

PROOF. The statement of the lemma follows from Exercise I.5.4 of [2]. □

LEMMA 5.2. (a) *The eigenvalues of $Q^{(n)}(z)$ are the zeros of $P_n(x; z)$.*

(b) *The eigenvalues of $\hat{Q}^{(n)}(z)$ are the zeros of $\hat{P}_n(x; z)$ and for each such eigenvalue $x$, the corresponding left eigenvector is given by*

$$(P_0(x; z), P_1(x; z), \ldots, P_{n-1}(x; z)).$$

Proof. We have already observed [after (43)] that 0 is an eigenvalue of $Q^{(n)}(z)$ if and only if it is a zero of $P_n(x; z)$. For the general case, let $I_n$ denote the identity matrix of dimension $n$. Write $Q^{(n)}$ for $Q^{(n)}(z)$ and similarly for $\widehat{Q}^{(n)}(z)$. The characteristic polynomial of $Q^{(1)}$ is

$$\det(xI_1 - Q^{(1)}) = x + \lambda + \mu_2(1 - z).$$

Because the $Q^{(n)}$ are tridiagonal, we have

$$\det(xI_2 - Q^{(2)}) = (x + \lambda + \mu_1 + \mu_2(1 - z))\det(xI_1 - Q^{(1)}) - \frac{\mu_1}{z}\lambda$$

and, for $n \geq 3$,

$$\det(xI_n - Q^{(n)}) = (x + \lambda + \mu_1 + \mu_2(1 - z))\det(xI_{n-1} - Q^{(n-1)})$$
$$- \frac{\mu_1}{z}\lambda\det(xI_{n-2} - Q^{(n-2)}).$$

Hence, we see that $(\mu_1/z)^n P_n(x; z)$ is the characteristic polynomial of $Q^{(n)}$, and thus, for each $n \geq 1$, the eigenvalues of $Q^{(n)}$ are the zeros of $P_n(x; z)$. This proves (a).

To show the first part of (b), observe that the characteristic polynomial of $\widehat{Q}^{(n)}$ satisfies

$$\det(xI_n - \widehat{Q}^{(n)}) = (x + \mu_1 + \mu_2(1 - z))\det(xI_{n-1} - Q^{(n-1)})$$
$$- \frac{\mu_1}{z}\lambda\det(xI_{n-2} - Q^{(n-2)})$$
$$= \det(xI_n - Q^{(n)}) - \lambda\det(xI_{n-1} - Q^{(n-1)})$$
$$= \left(\frac{\mu_1}{z}\right)^n\left(P_n(x; z) - \frac{\lambda z}{\mu_1}P_{n-1}(x; z)\right).$$

Hence, the eigenvalues of $\widehat{Q}^{(n)}$ are the zeros of $\widehat{P}_n(x; z)$.

To prove the second part of (b), it is readily checked that for each eigenvalue $\hat{x}$ of $\widehat{Q}^{(n)}$, for which $P_n(\hat{x}; z) = \lambda z P_{n-1}(\hat{x}; z)/\mu_1$, we have

$$(P_0(\hat{x}; z), P_1(\hat{x}; z), \ldots, P_{n-1}(\hat{x}; z))(\hat{x}I_n - \widehat{Q}^{(n)}) = \mathbf{0}. \qquad \square$$

Since $\widetilde{R}_m$ is positive (see Remark 3.1), an eigenvector $\mathbf{w} = (w_0, \mathbf{w}_1)$ of $R_m$ can be positive if and only if $\mathbf{w}_1$ is the Perron–Frobenius eigenvector of $\widetilde{R}_m$. By Theorem 2.5, $\mathbf{w}$ is an eigenvector of $\widehat{Q}^{(m+1)}(z)$ with eigenvalue zero and, because $\widehat{Q}^{(m+1)}(z)$ is an ML-matrix (see [16]), $\mathbf{w}$ can be positive if and only if zero is the *largest* eigenvalue of $\widehat{Q}^{(m+1)}(z)$. In Lemma 5.3, we shall show that there is exactly one $z \in (0, 1)$ such that the largest eigenvalue of $\widehat{Q}^{(m+1)}(z)$ is zero.

LEMMA 5.3.   *For $m \geq 1$ there exists a unique number $\hat{z}_{m+1}$ in the interval $(0,1)$ such that $\hat{x}_{m+1,m+1}(\hat{z}_{m+1}) = 0$.*

PROOF.   Consider the nonnegative matrix

$$(53) \qquad \Xi_{m+1}(z) \equiv \frac{z(\lambda + \mu_1 + \mu_2)I_{m+1} + z\hat{Q}^{(m+1)}(z)}{\lambda + \mu_1 + \mu_2}$$

and let $\xi_{m+1}(z)$ denote its largest eigenvalue. For $z \in (0,1)$, $\Xi_{m+1}(z)$ is a substochastic matrix, which is stochastic when $z = 1$. Thus Lemma 1.3.4 of [11] can be applied. Specifically, under the appropriate stability condition (23) or (24), the equation $z = \xi_{m+1}(z)$ has exactly one solution $\hat{z}_{m+1} \in (0,1)$. It is readily seen that $z = \xi_{m+1}(z)$ if and only if the maximum eigenvalue of $\hat{Q}^{(m+1)}(z)$ is equal to zero and the result follows.   $\square$

We have now proved the following theorem.

THEOREM 5.4.   *When $m$ is finite, the maximal eigenvalue of $R_m$ is given by the unique $\hat{z}_{m+1} \in (0,1)$ such that $\hat{x}_{m+1,m+1}(\hat{z}_{m+1}) = 0$. The corresponding eigenvector is strictly positive. The eigenvectors corresponding to any other nonzero eigenvalue of $R_m$ cannot be nonnegative.*

In view of (12), it is obvious that Theorem 5.4 determines the geometric decay rate of the level process we were looking for. The following corollary concerns the limiting behavior of this decay rate as $m$, the size of the first buffer, tends to infinity.

COROLLARY 5.5.   *Let $r_m$ be the Perron–Frobenius eigenvalue of $R_m$ for finite $m$.*
*If $\mu_1 \leq \mu_2$, then $r_1, r_2, \ldots$ strictly increases to $\eta$.*
*On the other hand, if $\mu_1 > \mu_2$, then $r_1, r_2, \ldots$ strictly increases to $\rho_2$.*

PROOF.   It was stated in Lemma 4.6 that $\{x_{n,n}(z)\}$ strictly increases to $\chi_1(z)$. To prove that $\{\hat{x}_{n,n}(z)\}$ also increases to $\chi_1(z)$, the interlacing property ensures that we need only to show $\hat{P}_n(\chi_1(z); z) > 0$ for $n \geq 1$.

For the case $z \geq \rho_1$ we have $\chi_1(z) = \chi(z)$, and from (48)–(51) it is easily checked by induction that $\hat{P}_n(\chi(z); z) = (1-z)(\lambda/\mu_1)^n > 0$.

For the case $z < \rho_1$, where $\chi_1(z) = \tau(z)$, first note that

$$P_n(\tau(z); z) - \sqrt{z\rho_1}P_{n-1}(\tau(z); z) > 0.$$

This can be shown easily by induction, using (41)–(43). Since we can write

$$\hat{P}_n(\tau(z); z) = P_n(\tau(z); z) - \sqrt{z\rho_1}P_{n-1}(\tau(z); z) + (\sqrt{z\rho_1} - z\rho_1)P_{n-1}(\tau(z); z),$$

$0 < z\rho_1 < 1$ and $P_{n-1}(\tau(z); z) > 0$ (see Lemma 4.7), we conclude that $\widehat{P}_n(\tau(z); z) > 0$.

Now, by Lemmas 5.3 and 4.5, the sequence $r_m = \{\hat{z}_{m+1}\}$ increases strictly to a $z^* \in (0, 1)$ which is the unique zero of $\chi_1(z)$ in the interval $(0, 1)$. Assume that $\mu_1 \leq \mu_2$. By Lemma 4.1, this can occur only when $\eta \leq \rho_1$. In this case $\chi_1(z) = \tau(z)$, which has a zero at $z = \eta$. Thus $r_m = \{\hat{z}_{m+1}\}$ increases strictly to $\eta$. On the other hand, when $\mu_1 > \mu_2$, Lemma 4.1 implies that $z > \rho_1$ and $\chi_1(z) = \chi(z)$, which has its zero at $z = \rho_2$. The sequence $r_m = \{\hat{z}_{m+1}\}$ then increases to $\rho_2$. $\square$

The above result shows that we must clearly distinguish between two possible regimes. These correspond with the different cases identified in Lemma 4.1. In the first regime, when $\mu_1 \leq \mu_2$, the first queue is the bottleneck and $\lim_{m\to\infty} \mathrm{sp}(R_m) = \eta$. In the second regime, when $\mu_1 > \mu_2$, the second queue is the bottleneck and $\lim_{m\to\infty} \mathrm{sp}(R_m) = \rho_2$. Note also that, in this second regime, Lemma 4.6 tells us that the limit of the sequence of the maximal eigenvalues of $R_m$ is different from the limit of the sequence of second-largest eigenvalues and so the limiting spectrum of $R_m$ has an isolated point.

We observed in Section 4 that the decay rate of the tandem Jackson network with infinite waiting room at the first queue is always $\rho_2$ irrespective of whether $\mu_1 \leq \mu_2$ or $\mu_1 > \mu_2$. We thus see that, when $\mu_1 \geq \mu_2$, the limiting decay rate of the finite truncations is indeed that of the infinite system. However, if $\mu_1 < \mu_2$, the limiting decay rate of the finite truncations is different from that of the infinite system. We have thus provided a counterexample to the idea that the decay rate of a QBD process with infinitely many phases can be derived by calculating the decay rates of finite truncations and then allowing the point at which truncation occurs to grow to infinity.

**6. Varying the decay rate.** An interesting question arises from the observations at the end of Section 4. By appropriately changing the transition intensities at level zero, in other words changing the entries in $\widetilde{Q}_1$, can we ensure that the stationary distribution decays at a rate that is given by any of the feasible values of $z$? In changing $\widetilde{Q}_1$ we have a great deal of freedom, so we might expect that the answer is yes. In fact it is. Below, we present two examples in which $\widetilde{Q}_1$ remains a tridiagonal matrix.

EXAMPLE 6.1. Suppose $\mu_1 > \mu_2$. We wish to have a decay rate $z$, satisfying the conditions in Theorem 4.9, which in this case means that $z \in [\rho_2, 1)$. By Lemmas 4.2 and 4.8 the vector $\mathbf{w}$ given in (34) will be positive and in $\ell^1$. We now replace each $\lambda$ in $\widetilde{Q}_1$ by a phase-dependent $\tilde{\lambda}_i$. Specifically, we define $\tilde{\lambda}_i$ recursively by

$$
\begin{aligned}
(54) \qquad & \tilde{\lambda}_0 = \mu_2 z, \\
& \tilde{\lambda}_i = \tilde{\lambda}_{i-1}\frac{w_{i-1}}{w_i} + \mu_2 z - \mu_1, \qquad i = 1, 2, \dots.
\end{aligned}
$$

The following proposition shows that this defines proper transition intensities.

PROPOSITION 6.1.   *The sequence $\{\tilde{\lambda}_i\}_{i=0}^{\infty}$ is strictly positive.*

PROOF.   Let $\phi(v)$ be the generating function of the sequence $w_1, w_2, \dots$. From (32) we find after some algebra that

$$\phi(v) = \frac{\mu_1(1 - vz)}{\lambda z v^2 - z(\lambda + \mu_1 + \mu_2(1 - z))v + \mu_1}.$$

Substituting $v = 1$ gives

(55)
$$\phi(1) = \sum_{i=0}^{\infty} w_i = \frac{\mu_1}{\mu_1 - \mu_2 z}.$$

Now, consider the sequence $y_0, y_1, \dots,$ with $y_i = \tilde{\lambda}_i w_i$. This sequence satisfies the recursion

$$y_i = y_{i-1} + (\mu_2 z - \mu_1) w_i, \qquad i = 1, 2, \dots,$$

with $y_0 = \mu_2 z > 0$. If $\mu_2 z \geq \mu_1$, then all $y_i$ (and hence $\tilde{\lambda}_i$) are obviously positive. On the other hand, if $\mu_2 z < \mu_1$, then $y_1, y_2, \dots$ is monotone decreasing, with

$$\lim_{i \to \infty} y_i = \mu_2 z + (\mu_2 z - \mu_1) \sum_{i=1}^{\infty} w_i = 0,$$

which shows that all $\tilde{\lambda}_i$ are positive in this case as well.   □

The recursion (54) ensures that $\mathbf{w}$ is a $(\mu_2 z)^{-1}$-invariant measure of $R_\infty Q_2$. Moreover, $\mathbf{w}$ satisfies $\mathbf{w}\tilde{Q}_1 = -\mu_2 z \mathbf{w}$. Hence $\mathbf{w}(\tilde{Q}_1 + R_\infty Q_2) = \mathbf{0}$, so that by (5) and Theorem 2.2 it follows that the stationary distribution $\boldsymbol{\pi} = (\boldsymbol{\pi}_0, \boldsymbol{\pi}_1, \dots)$ of $(Y_t, J_t)$ is given by

$$\boldsymbol{\pi}_n = c\mathbf{w}R_\infty^n = z^n c\mathbf{w}, \qquad n \geq 0,$$

for some normalizing constant $c$. Thus, it is clear that $z$ indeed is the decay rate in this model.

This example has demonstrated the counterintuitive result that, by changing the arrival intensity to the first queue when the second is empty, such that it becomes dependent of the number of customers in the first queue, we can produce any decay rate in the range $[\rho_2, 1)$.

EXAMPLE 6.2.   Suppose $\mu_1 < \mu_2$. We wish to have a decay rate $z$, with $z \in [\eta, \rho_2]$. Again, the vector $\mathbf{w}$ given in (34) is positive and in $\ell^1$. This time we leave the arrival rate unchanged, but introduce an extra transition rate $\nu_i$ from state $(0, i)$ to $(0, i-1)$. This corresponds to removing customers from the first queue, without introducing them to the second queue. The values $\nu_i$ are recursively defined as

$$\nu_1 = \frac{(\lambda - \mu_2 z) w_0}{w_1},$$

$$\nu_{i+1} = \frac{(\nu_i + \lambda + \mu_1 - \mu_2 z) w_i - \lambda w_{i-1}}{w_{i+1}}.$$

PROPOSITION 6.2. *The sequence $\{\nu_i\}_{i=0}^{\infty}$ positive.*

PROOF. The proof is similar to the proof of Proposition 6.1. First, we claim that

(56) $$(\lambda + \mu_1 - \mu_2 z) w_i < \lambda w_{i-1}, \qquad i = 1, 2, \ldots.$$

To see this, consider the sequence of polynomials $\{\Theta_n\}$, defined by $\Theta_n(x) = (\lambda + \mu_1 - \mu_2 z) P_n(x) - \lambda P_{n-1}(x)$, $n \geq 1$, with the polynomials $\{P_n\}$ given in (41)–(43). Imitating the proof of Lemma 4.5 for $\Theta_n$ instead of $\widehat{P}_n$, we find that the zeros of $\{\Theta_n\}$ interlace, that the largest zero of $\Theta_n$ is larger than the largest zero of $P_n$ and that the second largest zero of $\Theta_n$ is smaller than the largest zero of $P_n$. Now, for $z \in [\eta, \rho_2]$, the largest zero of $P_n$ is less than or equal to 0. Hence, $\Theta_n$ can have at most one zero greater than 0. It is easily verified that the largest zero of $\Theta_1$ is given by $(\lambda - \mu_2 z)(\mu_2 z^2 - (\lambda + \mu_1 + \mu_2) z + \mu_1)/(z(\lambda - \mu_2 z + \mu_1))$, which is strictly positive for all $0 < z \leq \rho_2$. Hence, all $\Theta_n$ have exactly one strictly positive zero. Thus, because the leading coefficient of $\Theta_n$ is positive, $\Theta_n(0)$ must be strictly negative, which is equivalent to (56).

Second, let $y_i = \lambda_i w_i$, $i = 1, 2, \ldots$. We have, for all $i = 2, 3, \ldots$,

$$y_i = y_{i-1} + (\lambda + \mu_1 - \mu_2 z) w_i - \lambda w_{i-1},$$

where $y_1 = \lambda - \mu_2 z > 0$. Thus, using (56), $y_1, y_2, \ldots$ is a strictly decreasing sequence with limit

$$\lambda - \mu_2 z + (\lambda + \mu_1 - \mu_2 z) \sum_{i=1}^{\infty} w_i - \lambda \sum_{i=0}^{\infty} w_i = 0,$$

where we have again used (55). This shows that all $\nu_i$ are positive. □

As above, the recursion ensures that $\mathbf{w}(\widetilde{Q}_1 + R_\infty Q_2) = \mathbf{0}$, so that the stationary distribution of $(Y_t, J_t)$ is given by

$$\boldsymbol{\pi}_n = c \mathbf{w} R_\infty^n = c z^n \mathbf{w}, \qquad n \geq 0,$$

for some normalizing constant $c$, from which it is clear that $z$ is the decay rate in this model.

Thus, by allowing customers at the first queue to be removed at specified rates when the second queue is empty, we have been able to produce any decay rate in $[\eta, \rho_2]$. Note that it is not possible using this scheme to produce a decay rate greater than $\rho_2$. However, we can do so using a scheme such as that in Example 6.1, which is also applicable here, since the proof of Proposition 6.1 did not use the fact that $\mu_1 > \mu_2$.

In the examples given above, $\widetilde{Q}_1$ was constructed such that the corresponding $\boldsymbol{\pi}_0$ is exactly equal to some $z^{-1}$-invariant measure of $R_\infty$. As a consequence, the stationary distribution of $(Y_t, J_t)$ has a product form. However, it is also possible to construct $\widetilde{Q}_1$, so that $\boldsymbol{\pi}_0$ is a finite linear combination of $z^{-1}$-invariant measures of $R_\infty$. In that case the stationary distribution does not have a product form. The decay rate is then given by the largest value of $z$ with corresponding $z^{-1}$-invariant measure in the linear combination.

REMARK 6.3.   When $\mu_1 \leq \mu_2$, the minimal attainable decay rate cannot be less than $\eta$ and when $\mu_1 > \mu_2$ the minimal attainable decay rate cannot be less than $\rho_2$. This follows because $\eta$ and $\rho_2$ respectively are the smallest values of $z$ for which a $z^{-1}$-invariant measure exists.

The maximal attainable decay rate is produced in a different way. Clearly, when $\mu_1 > \mu_2$, any decay rate in $[\rho_2, 1)$ can be produced. However, when $\mu_1 < \mu_2$, it is not immediately clear whether the matrix $R_\infty$ has a value $z \in [\mu_1/\mu_2, 1)$ with a corresponding $z^{-1}$-invariant measure $\mathbf{w}$ that is not in $\ell^1$. If such a measure did exist, the behavior at level 0 would be such that the first queue is unstable, while the second remains stable and has decay rate $z$. A physical argument tells us that this is possible only when $z = \mu_1/\mu_2$: if the first queue is unstable, then the second queue behaves like a standard $M/M/1$ queueing system with arrival rate $\mu_1$ and service rate $\mu_2$. This implies that its decay rate could never be larger than $\mu_1/\mu_2$.

**7. Hitting probabilities on high levels: general QBD processes.**   In various applications one is interested in *hitting* or *exit* probabilities of the level process. In this and the following section, we shall consider the decay rate of these probabilities, first in the context of a general QBD with possibly infinitely many phases and then in the context of the $M/M/1$ tandem.

For $a < b$, define $T_a^b$ to be the first time that either level $a$ or level $b$ is hit. Also, let $\mathbb{P}_{ki}$ denote the probability measure under which the QBD process starts in $(k, i)$. For $k \geq 0$, we are interested in the decay rate as $K \to \infty$ of the first exit probabilities

$$P_k^K(i, j) := \mathbb{P}_{ki}(J_{T_0^K} = j, Y_{T_0^K} = K),$$

which we collect into a matrix $P_k^K$. Define the matrix $H_k$ to be equal to $P_k^{k+1}$. Thus, $H_0$ is the 0 matrix and it is not difficult to see that, for $k \geq 1$, $H_1, H_2, \ldots$ satisfy the recursion

$$Q_0 + Q_1 H_k + Q_2 H_{k-1} H_k = 0,$$

and the matrix $P_k^K$ is given by

(57) $$P_k^K = H_k H_{k+1} \cdots H_{K-1}.$$

The following result is essentially a restatement of Lemma 8.2.1 of [7].

LEMMA 7.1. *The sequence of matrices, $H_1, H_2, \ldots$ increases elementwise to the matrix $H$ which is the minimal nonnegative solution to the matrix equation*

(58) $$Q_0 + Q_1 H + Q_2 H^2 = 0.$$

For the case when $m$ is finite, it was shown in [9], Lemma 3.1, that, when Condition 2.1 holds, either $H$ is primitive or, by a suitable permutation of the states, it can be written in the form

(59) $$H = \begin{bmatrix} L_1 & 0 \\ L_{.1} & L_. \end{bmatrix},$$

where $L_1$ is primitive and $L_.$ is lower triangular with its diagonal entries equal to zero. A similar result can be established even when $m$ is infinite. (At the time of writing, this result, due to Latouche and Taylor, is unpublished. An explanation can be obtained from Peter Taylor at p.taylor@ms.unimelb.edu.au.) Thus, when $m \leq \infty$, $H$ has the decomposition (59), where $L_.$ is lower triangular and $L_1$ is irreducible and aperiodic. It follows from (14) that the convergence norm $c$ of $L_1$ is well defined and given by

$$c = \lim_{n \to \infty} (L_1^n(i, j))^{1/n}.$$

Let $\Sigma^*$ be the set of indices corresponding to $L_1$ and partition the matrices $H_k$, conformally with our partition of the matrix $H$, so that

(60) $$H_k = \begin{bmatrix} L_1^{(k)} & 0 \\ L_{.1}^{(k)} & L_.^{(k)} \end{bmatrix}.$$

The decay behavior of the hitting probabilities is described in Theorem 7.3. However, first we need a lemma.

LEMMA 7.2. *For any phase $i$, there exist numbers $k^*$ and $N^*$ such that, for $k > k^*$ and $N > \max(k, N^*)$, there is a $\nu \in \Sigma^*$ with*

(61) $$P_k^N(i, \nu) > 0.$$

PROOF. First note that, for any given $N$ and $k$, say $N_0$ and $k_0$, there may not be a $\nu_0 \in \Sigma^*$ such that $P_{k_0}^{N_0}(i, \nu_0) > 0$. Taking into account the irreducibility of the doubly infinite process with generator (2), this could be because every path of positive probability from state $(k_0, i)$ to states of the form $(N_0, \nu_0)$ with $\nu_0 \in \Sigma^*$ does one of the following:

1. passes through a state of the form $(N_0, m)$ with $m \notin \Sigma^*$ [by the decomposition (60), this can occur only if $i \notin \Sigma^*$];
2. goes through level zero;
3. does both 1 and 2.

Consider a path from $(k, i)$ to a state $(N_0, \nu_0)$ with $N_0 > k$ and $\nu_0 \in \Sigma^*$ of the form described in 1 above. Let $N^*$ be the highest level it reaches. Any path from $(N_0, \nu_0)$ to level $N^* + 1$ must hit level $N^* + 1$ in a state $(N^* + 1, \nu)$ with $\nu \in \Sigma^*$ and, by irreducibility of the process with generator (2), there must be such a path. Concatenating these two paths, we have constructed a path of positive probability from $(k, i)$ to $(N^* + 1, \nu)$ which first hits level $N^* + 1$ in phase $\nu \in \Sigma^*$. If this path does not pass through level 0, then we have constructed a path as desired. If it does go through level 0 we modify it as described below.

If, after performing the modification described above, all paths from state $(k_0, i)$ to states of the form $(N_0, \nu_0)$ with $\nu_0 \in \Sigma^*$ pass through level 0 choose one such path let $-\tilde{k}$ be the lowest level reached by the path and put $k^* = k_0 + \tilde{k}$. Then if $k > k^*$, there is a path of positive probability from $(k, i)$ to a state of the form $(N_0 + k - k_0, \nu_0)$ with $\nu_0 \in \Sigma^*$ which does not pass through level 0. This shows that, for all phases $i$, we can choose $k^*$ such that when $k > k^*$ there exists a path of positive probability which does not pass through level 0 from state $(k, i)$ to any level $N > \max(k, N_0)$. The lemma is thus proved. □

THEOREM 7.3. *Consider an irreducible QBD process with a finite or infinite phase space, satisfying Condition* 2.1 *and* (11). *Then for* $i \in \{0, 1, \dots\}$ *and* $j \in \Sigma^*$ *there exists* $k^*$ *such that, for* $k > k^*$,

$$\lim_{K \to \infty} \frac{\log P_k^K(i, j)}{K} = \log(c), \tag{62}$$

*where* $c$ *is the convergence norm of* $L_1$. *For* $i, j \notin \Sigma^*$, *there exists* $K^*$ *such that, for* $K > K^*$,

$$P_k^K(i, j) = 0. \tag{63}$$

*For* $i \in \Sigma^*$, $j \notin \Sigma^*$ *and all* $k < K$,

$$P_k^K(i, j) = 0. \tag{64}$$

PROOF.   For the case where $i, j \in \Sigma^*$, we have

$$\frac{\log(P_k^K(i,j))}{K} = \frac{\log(L_1^{(k)} \cdots L_1^{(K-1)}(i,j))}{K}$$

$$\leq \frac{\log((L_1)^{K-k}(i,j))}{K}$$

$$= \frac{K-k}{K} \log((L_1)^{K-k}(i,j)^{1/(K-k)})$$

so that, letting $K \to \infty$, we find by (14) that

(65) $$\limsup_{K \to \infty} \frac{\log(P_k^K(i,j))}{K} \leq \log(c).$$

To show the opposite, choose $k^*$, $N^*$ and $\nu$ so that (61) is satisfied. Then, for $k > k^*$, $N > \max(k, N^*)$ and $K > N$, we have

(66) $$P_k^K(i,j) = \sum_l P_k^N(i,l) L_1^{(N)} \cdots L_1^{(K-1)}(l,j)$$

$$\geq P_k^N(i,\nu)(L_1^{(N)})^{K-N}(\nu,j).$$

Now we have

$$\frac{\log(P_k^K(i,j))}{K} \geq \frac{\log(P_k^N(i,\nu))}{K} + \frac{\log((L_1^{(N)})^{K-N}(\nu,j))}{K}$$

$$= \frac{\log(P_k^N(i,\nu))}{K} + \frac{K-N}{K} \log((L_1^{(N)})^{K-N}(\nu,j)^{1/(K-N)}),$$

so that, letting $K \to \infty$, we find that

$$\liminf_{K \to \infty} \frac{\log(P_k^K(i,j))}{K} \geq \log(c_N),$$

where $c_N$ is the convergence norm of $L_1^{(N)}$. Since this holds for all $N$ and, by Lemma 2.3, $c_N \to c$ as $N \to \infty$, we see that

(67) $$\liminf_{K \to \infty} \frac{\log(P_k^K(i,j))}{K} \geq \log(c),$$

which, together with (65), gives the result in this case.

When $i \notin \Sigma^*$ and $j \in \Sigma^*$, we can still use Lemma 7.2 to choose $k^*$ and $N^*$ so that, when $N > N^*$ and $k > k^*$, there exists a $\nu \in \Sigma^*$ such that (61) is satisfied, and argue as above from (66) that (67) is satisfied.

To get the analogue of (65), observe that, by (57) and (60), we must be able to write

$$P_k^K(i,j) = \sum_{r=k}^{K-1} L^{(k)} \cdots L^{(r-1)} L_{\cdot 1}^{(r)} L_1^{(r+1)} \cdots L_1^{(K-1)}(i,j)$$

$$\leq \sum_{r=k}^{K-1} L_{.}^{r-k} L_{.1} L_1^{K-1-r}(i,j).$$

Now, because $L_{.}$ is lower triangular, there is a positive integer $s^*$ such that $L_{.}^s(i, \nu) = 0$ for all $s > s^*$ and $\nu \notin \Sigma^*$. Thus, for $K > k + s^* + 1$,

$$P_k^K(i,j) \leq \sum_{s=0}^{s^*} L_{.}^s L_{.1} L_1^{K-1-k-s}(i,j)$$

$$= (D L_1^{K-1-k-s^*})(i,j),$$

where

$$D = \sum_{s=0}^{s^*} L_{.}^s L_{.1} L_1^{s^*-s}.$$

Consider the Markov chain with transition matrix $H$. Let $\tau(j)$ be the first time greater than or equal to $s^* + 1$ that the chain visits state $j \in \Sigma^*$ and let $f^{(n)}(i,j)$ be the probability that $\tau(j) = n$, conditional on the chain starting in state $i \notin \Sigma^*$. Then it follows easily that

$$(D L_1^{K-1-k-s^*})(i,j) = \sum_{n=s^*+1}^{K-k} f^{(n)}(i,j) L_1^{K-k-n}(j,j)$$

and that

$$(68) \qquad DL(i,j;z) = \frac{F(i,j;z)L(j,j;z)}{z^{s^*+1}},$$

where

$$DL(i,j;z) = \sum_{n=0}^{\infty} (D L_1^n)(i,j) z^n,$$

$$F(i,j;z) = \sum_{n=s^*+1}^{\infty} f^{(n)}(i,j) z^n$$

and

$$L(j,j;z) = \sum_{n=0}^{\infty} L_1^n(j,j) z^n.$$

It is clear that $f^{(n)}(i,j) \leq (D L_1^n)(i,j)$ and so the convergence radius of the power series $F(i,j;z)$ is greater than or equal to the convergence radius of the power series $DL(i,j;z)$. Therefore, by (68), the convergence radii of the series $DL(i,j;z)$ and $L(j,j;z)$ are the same. Thus we have

$$\lim_{k \to \infty} \frac{\log((D L_1^{K-1-k-s^*})(i,j))}{K} = \log(c)$$

and so
$$\limsup_{K \to \infty} \frac{\log(P_k^K(i,j))}{K} \le \log(c).$$
Thus, the result is proved for $i \notin \Sigma^*$ and $j \in \Sigma^*$.

When $i, j \notin \Sigma^*$, using the same definition of $s^*$ as above, it follows that, for all $K > k + s^*$,

$$(69) \qquad P_k^K(i,j) = 0,$$

while we immediately have $P_k^K(i,j) = 0$ when $i \in \Sigma^*$ and $j \notin \Sigma^*$. This proves the second part of the theorem. □

To finish off this section, we present some results for the matrix $H_m$ that are analogous to Theorem 2.5 and Corollary 2.6 for the matrix $R_m$. This will allow us to conclude that the eigenvalues of $R_m$ and $H_m$ coincide when $m < \infty$; see Corollary 7.6.

As for (21) it is easy to see that, when the column vector $\mathbf{v}$ and scalar $z$ satisfy $H_m \mathbf{v} = z\mathbf{v}$, then

$$(70) \qquad (Q_0 + zQ_1 + z^2 Q_2)\mathbf{v} = \mathbf{0}.$$

Again, under certain conditions on $\mathbf{v}$ and $z$, the converse is true, irrespective of whether $m < \infty$ or $m = \infty$. This is shown in the following theorem, which is basically Theorem 5.3 of [14].

THEOREM 7.4. *Consider a continuous-time QBD process with generator of the form* (1). *Then, if the complex variable $z$ and the vector $\mathbf{v} = \{v_k\}$ are such that $|z| < 1$ and $\sum_k |v_k| q_k < \infty$, then* (70) *implies that*
$$H_m \mathbf{v} = z\mathbf{v}.$$

Following essentially the same proof as for Corollary 2.6 we have a characterization for the case $m < \infty$.

COROLLARY 7.5. *For an irreducible QBD process with a finite phase space, satisfying* (11), *the eigenvalues of $H_m$ are all the zeros of the polynomial*

$$(71) \qquad \det(Q_0 + zQ_1 + z^2 Q_2)$$

*that lie strictly within the unit circle.*

The following result is now immediate from Corollaries 2.6 and 7.5.

COROLLARY 7.6. *For an irreducible QBD process with a finite phase space, satisfying* (11), *the eigenvalues of $H_m$ and $R_m$ coincide.*

In particular, when $m < \infty$ the hitting probabilities $P_k^K(i,j)$ have the same geometric rate of decay as the stationary probabilities $\pi_{Ki}$ in (12).

**8. Hitting probabilities on high levels: the tandem network.** Assume that the tandem queue starts in state $(1, i)$ with $i \geq 1$, that is, with one customer in the second queue and $i \geq 1$ customers in the first queue. It is possible that the process can first hit level 2 before level 0 with any number $j \geq i-1$ customers in the first queue, that is, in any state $(2, j)$ with $j \geq i-1$. If the queue starts in state $(1, 0)$, then it can first hit level 2 before level 0 with any number of customers in the first queue. A consequence of this is that the matrix $H_m$ for this QBD process is irreducible. It then follows from Theorem 7.3 that

$$\lim_{K \to \infty} \frac{\log P_1^K(i,j)}{K} = \log(c), \tag{72}$$

where $c$ is the convergence norm of $H_m$. Thus, to calculate the decay rate of the hitting probabilities, we need to calculate the convergence norm of $H_m$.

To do this for the case $m = \infty$, we could follow a line of reasoning similar to that we used in Section 4 based upon Theorem 7.4 instead of Theorem 2.5. Thus, we would calculate conditions for a solution $\mathbf{v}$ to $Q(z)\mathbf{v} = \mathbf{0}$ to be both positive and in $\ell^1$. However, unlike the $z^{-1}$-invariant measure of $R_\infty$, which affects the decay rate of the stationary distribution, the $z^{-1}$-invariant vector of $H_\infty$ has no effect on the decay rate of the hitting probabilities. We thus choose to calculate the decay rate of the hitting probabilities in a more efficient way.

By Corollary 7.6, for finite $m$, the eigenvalues of $R_m$ and $H_m$ coincide. Thus Theorem 5.4 and Corollary 5.5 apply to $H_m$ as well as to $R_m$. In particular, we have the following theorem.

THEOREM 8.1. *When $m$ is finite, the following results hold:*

1. *The maximal eigenvalue $h_m$ of $H_m$ is given by the unique $\hat{z}_{m+1} \in (0, 1)$ such that $\hat{x}_{m+1,m+1}(\hat{z}_{m+1}) = 0$.*
2. *If $\mu_1 \leq \mu_2$, then $h_1, h_2, \dots$ strictly increases to $\eta$. On the other hand, if $\mu_1 > \mu_2$, then $h_1, h_2, \dots$ strictly increases to $\rho_2$.*

By Lemma 2.3, it follows that the convergence norm of $H_\infty$ is then equal to $\eta$ if $\mu_1 \leq \mu_2$ and $\rho_2$ if $\mu_1 > \mu_2$. Together with Theorem 7.3, this gives us the following theorem.

THEOREM 8.2.   1. *When $m$ is finite,*

$$\lim_{K \to \infty} \frac{\log P_1^K(i,j)}{K} = \log(h_m). \tag{73}$$

2. *When $m$ is infinite:*

(a) *when $\mu_1 \leq \mu_2$,*

$$(74) \qquad \lim_{K \to \infty} \frac{\log P_1^K(i,j)}{K} = \log(\eta);$$

(b) *when $\mu_1 > \mu_2$,*

$$(75) \qquad \lim_{K \to \infty} \frac{\log P_1^K(i,j)}{K} = \log(\rho_2).$$

The decay rate of the hitting probabilities in the case of infinite $m$ is thus the same as the decay rate of the stationary number in the second queue when $\mu_1 \geq \mu_2$, but it is not the same when $\mu_1 < \mu_2$. This is an interesting property of the tandem Jackson network, which we believe was not known previously.

There are two further interesting questions about the decay rate of the hitting probabilities that we have not addressed above. The first question involves the decay rate of $\sum_j P_k^K(i,j)$ as $K \to \infty$ in the case $m = \infty$. It follows immediately from Theorem 8.2 that this decay rate is larger than $\eta$ and $\rho_2$ when $\mu_1 \leq \mu_2$ and $\mu_1 > \mu_2$, respectively. We conjecture that it is equal to these values, although we currently have no proof of this.

The second question involves the decay rate of the hitting probabilities on level $K$ if the process starts in level 1 according to some distribution $\mathbf{x}_1$. Of particular interest is the situation when $\mathbf{x}_1$ is the stationary distribution $\boldsymbol{\pi}_1$ at level 1. These hitting probabilities are given by the components of

$$(76) \qquad \mathbf{x}_1 P_1^K.$$

When $m < \infty$, the decay rates of these probabilities are easily seen to be the same as the decay rates of $P_1^K(i,j)$, given by Theorem 8.2. However, when $m = \infty$, this need not necessarily hold. Indeed, our experience with the matrix $R_\infty$ would lead us to believe that we could achieve any decay rate in $[\eta, \mu_1/\mu_2)$ if $\mu_1 \leq \mu_2$ and any decay rate in $[\rho_2, 1)$ if $\mu_1 > \mu_2$. However, since we have no theorem analogous to Theorem 7.4 that can inform us about the $z^{-1}$-invariant measures, rather than the $z^{-1}$-invariant vectors, of $H_\infty$ we do not currently see how this problem can be approached.

**Acknowledgments.** Authors thank two anonymous referees for many helpful comments on an earlier version of this paper.

D. P. KROESE
DEPARTMENT OF MATHEMATICS
UNIVERSITY OF QUEENSLAND
QUEENSLAND 4072
AUSTRALIA

W. R. W. SCHEINHARDT
DEPARTMENT OF APPLIED MATHEMATICS
UNIVERSITY OF TWENTE
P.O. BOX 217
7500 AE ENSCHEDE
THE NETHERLANDS
AND
CENTRE FOR MATHEMATICS
   AND COMPUTER SCIENCE (CWI)
P.O. BOX 94079
1090 GB AMSTERDAM
THE NETHERLANDS

P. G. Taylor
Department of Mathematics
    and Statistics
University of Melbourne
Victoria 3010
Australia
e-mail: p.taylor@ms.unimelb.edu.au